\documentclass[12pt]{article}
\usepackage{amsmath}
\usepackage{amssymb}
\usepackage{amsbsy}
\usepackage{fancyhdr}
\usepackage{amsfonts}
\newtheorem{lem}{Lemma}[section]%
\newtheorem{theorem}[lem]{Theorem}%
\newtheorem{defi}[lem]{Definition}%
\newtheorem{cor}[lem]{Corollary}%
\newtheorem{prop}[lem]{Proposition}%

\def\a{\alpha}    
    \def\k{\kappa}

\def\G{\Gamma}

 \def\lg{\langle} \def\rg{\rangle}
\def\nd{\mathrel{\bigm|\kern-.7em/}}

\def\f{\noindent}

\def\Aut{\hbox{\rm Aut\,}}

\def\Cay{\hbox{\rm Cay }}

\def\demo{{\bf Proof}\hskip10pt}

\def\mz{{\mathbb Z}}
\def\qed{\hskip10pt $\Box$\vspace{3mm}}

\begin{document}
\title{On $g$-Extra Connectivity of Hypercube-like Networks
\author
{Jin-Xin Zhou\\
{\small Mathematics, Beijing Jiaotong University, Beijing 100044, P.R. China}\\
{\small\em{Email}:\ jxzhou@bjtu.edu.cn}}}

\date{}
\maketitle

\begin{abstract}

Given a connected graph $G$ and a non-negative integer $g$, the {\em $g$-extra connectivity} $\k_g(G)$ of $G$ is the minimum
cardinality of a set of vertices in $G$, if it exists, whose deletion disconnects $G$ and leaves each remaining component with
more than $g$ vertices. This paper focuses on the $g$-extra connectivity
of hypercube-like networks (HL-networks for short) which includes
numerous well-known topologies, such as hypercubes, twisted cubes, crossed cubes and M\"obius
cubes. All the known results suggest the equality $\k_g(X_n)=f_n(g)$ holds,
where $X_n$ is an $n$-dimensional HL-network, $f_n(g)=n(g+1)-\frac{g(g+3)}{2}$, $n\geq 5$ and
$0\leq g\leq n-3$? Some authors also attempted to prove this equality in general.
In this paper, we construct a subfamily of an $n$-dimensional HL-network with $g$-extra connectivity greater than
$f_n(g)$ which implies that the above equality does not hold in general.
We also prove that for $n\geq 5$ and
$0\leq g\leq n-3$, $\k_g(X_n)\geq f_n(g)$ always holds. This enables us to give a sufficient condition for the
equality $\k_g(X_n)=f_n(g)$, which is then used to determine the $g$-extra connectivity of HL-networks for some small
$g$ or the $g$-extra connectivity of some particular subfamily of HL-networks.
As a result, a short proof for the main results in [Journal of Computer and System Sciences 79 (2013) 669--688].



\bigskip

\end{abstract}

\section{Introduction}
With the rapid development of VLSI technology and software technology,
a multiprocessor system may contain hundreds
or even thousands of nodes. With the continuous increase
in the size of multiprocessor systems, the
complexity of a system can adversely affect its
fault tolerance and reliability. To the design and
maintenance purpose of multiprocessor systems, appropriate
measures of reliability should be found.

In a multiprocessor system, processors are connected based on a specific interconnection network. An interconnection
network is usually represented by a graph in which vertices represent processors and edges represent links between processors.
The traditional connectivity is an important factor for measuring the reliability of an interconnection network, which can correctly
reflects the fault tolerance of systems with few processor. However, it always underestimates
the resilience of large networks. The discrepancy incurred is because events
whose occurrence would disrupt a large network after a few processor failures are
highly unlikely, therefore, the disruption envisaged occurs in a worst case scenario. Motivated
by the shortcomings of the traditional connectivity, Harary~\cite{Harary} introduced the
concept of conditional connectivity.

Let $G$ be a connected undirected graph, and $\mathcal{P}$ a graph-theoretic property. Harary~\cite{Harary} defined the conditional connectivity
$\k(G;\mathcal{P})$ as the minimum cardinality
of a set of vertices, if any, whose deletion disconnects $G$ and
every remaining component has property $\mathcal{P}$.
Subsequently, F\'abrega and Fiol~\cite{FF} investigated the following kind of conditional connectivity.
A subset of vertices $S$ is said to be a {\em cutset} if $G-S$ is not connected. A cutset $S$ is called an $R_g$-cutset, where $g$ is
a non-negative integer, if every component of $G-S$ has at least $g+1$ vertices. If $G$ has at least one $R_g$-cutset, the {\em $g$-extra connectivity}
of $G$, denoted by $\k_g(G)$, is then defined as the minimum cardinality over all $R_g$-cutsets of $G$. In other
words, $\k_g(G)=\k(G;\mathcal{P}_g )$, where $\mathcal{P}_g$ denotes that every remaining component has more than $g$ vertices.

Obviously, $\k_0(G)=\k(G)$ for any
connected graph $G$ that is not a complete graph. Therefore, the $g$-extra
connectivity can
be regarded as a general form of the classical connectivity
that provides
measures that are more accurate for reliability and fault tolerance
for large-scale parallel processing systems. Regarding
the computational complexity of the problem,
based on thorough research, no polynomial-time algorithm
has been presented to compute $\k_g$ for a general
graph~\cite{CTH2013}; nor has there been any tight upper bound for
$\k_g$~\cite{Esfahanian0}. The problem of determining the $g$-extra
connectivity of numerous networks
has received a great deal of attention in recent years.
For more results regarding $g$-extra connectivity, see, for example, \cite{Boesch,CH,CTH2013,Esfahanian0,Esfahanian,FF1994,FF,Latifi,Wan-Z,YM,Zhu2008,Zhu}.
It is worthwhile to mention that different types of generalized connectivity such as
$g$-extra connectivity have many applications. One of them is the conditional diagnosability, which was firstly proposed by Lai
et al.~\cite{LTCH}. For extensive study, the readers may also refer to \cite{CH2010,FH,XTH}.

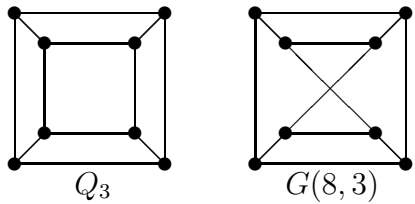
\begin{figure}[ht]
\begin{center}
\unitlength 4mm
\begin{picture}(15,6)

\put(1, 1){\circle*{0.4}} \put(3, 0){$Q_3$}

\put(1,1){\line(1, 0){5}} \put(1,1){\line(0, 1){5}}

\put(6, 1){\circle*{0.4}}

\put(6,1){\line(0, 1){5}}\put(1,6){\line(1, 0){5}}

\put(1, 6){\circle*{0.4}} 

\put(6, 6){\circle*{0.4}}


\put(2, 2){\circle*{0.4}} 

\put(2,2){\line(1, 0){3}} \put(2,2){\line(0, 1){3}}

\put(2,2){\line(-1, -1){1}}\put(5,2){\line(1, -1){1}}

\put(5, 2){\circle*{0.4}}

\put(5,2){\line(0, 1){3}}\put(2,5){\line(1, 0){3}}

\put(2, 5){\circle*{0.4}} 

\put(2,5){\line(-1, 1){1}}

\put(5, 5){\circle*{0.4}}

\put(5,5){\line(1, 1){1}}


\put(9, 1){\circle*{0.4}} \put(10, 0){$G(8,3)$}

\put(9,1){\line(1, 0){5}} \put(9,1){\line(0, 1){5}}

\put(14, 1){\circle*{0.4}}

\put(14,1){\line(0, 1){5}}\put(9,6){\line(1, 0){5}}

\put(9, 6){\circle*{0.4}} 

\put(14, 6){\circle*{0.4}}


\put(10, 2){\circle*{0.4}} 

\put(10,2){\line(1, 0){3}} \put(10,2){\line(1, 1){3}}

\put(10,2){\line(-1, -1){1}}\put(13,2){\line(1, -1){1}}

\put(13, 2){\circle*{0.4}}

\put(13,2){\line(-1, 1){3}}\put(10,5){\line(1, 0){3}}

\put(10, 5){\circle*{0.4}} 

\put(10,5){\line(-1, 1){1}}

\put(13, 5){\circle*{0.4}}

\put(13,5){\line(1, 1){1}}
\end{picture}
\end{center}\vspace{-.2cm}
\caption{Two $3$-dimensional HL-networks} \label{Fig-1}
\end{figure}

The {\em hypercube-like networks} ({\em HL-networks} for short) are defined recursively as follows:
$$\mathbb{L}_0=\{K_1\}\ {\rm and }\ \mathbb{L}_n=\{G_0\oplus G_1\ |\ G_0,G_1\in HL_{n-1}\},$$
where the symbol ``$\oplus$" represents the {\em perfect matching operation} that connects $G_0$ and $G_1$ using some perfect matching, denoted
by $PM(G)$. It is
obvious that $\mathbb{L}_1=\{K_2\}$, $\mathbb{L}_2=\{C_4\}$, and $\mathbb{L}_3=\{Q_3, G(8, 4)\}$, where $C_4$ is a cycle of length $4$, and
$Q_3$ and $G(8,4)$ are depicted as Figure~(\ref{Fig-1}).
(Some authors also use the term BC-networks instead~\cite{FH}.
In this paper, we follow~\cite{VRS} to use the term HL-networks.)
Numerous well-known topologies, such as hypercubes~\cite{SS}, crossed cubes~\cite{Efe}, M\"obius
cubes~\cite{cull}, twisted cubes~\cite{HKS}, varietal cubes~\cite{CC} etc. belong to the class of HL-networks.

For a positive integer $n$, let $f_n(g)=n(g+1)-\frac{1}{2}g(g+3)$ be a function of $g$.
From \cite{YM}, we know that $\k_g(Q_n)=f_n(g)$ for $n\geq 4$ and $0\leq g\leq n-3$.
Let $G$ be an $n$-dimensional HL-network. Fan et al. proved that $\k_0(G)=f_n(0)=n$ (see \cite{FH}),
and $\k_1(G)=f_n(1)=2n-2$ for $n\geq 3$ (see \cite{Zhu2008}). Xu et al.~\cite{XZX} proved that
$\k_2(G)=f_n(2)=3n-5$ for $n\geq 8$ and $\k_2(G)\geq 3n-5$ for $5\leq n\leq 7$. Recently, Chang and Hsieh~\cite{CH}
improved Xu's result by showing that $\k_2(G)=f_n(2)=3n-5$ for $n\geq 5$, and they also obtained that
$\k_3(G)=f_n(3)=4n-9$ for $n\geq 6$. The facts listed above provide a strong motivation for
studying the following problem. \medskip

\f{\bf Problem~A}\ For any $G\in \mathbb{L}_n$, does $\k_g(G)=f_n(g)$ hold for $0\leq g\leq n-3$?\medskip

In this paper, by analyzing the structure of
an $n$-dimensional HL-network $G$ with at most $f_n(g)$ faulty vertices, where $0\leq g\leq n-3$ and $n\geq 5$,
we give a lower bound on $\k_g(G)$, namely, $\k_g(G)\geq f_n(g)$. Furthermore, we give a sufficient
condition for the equality $\k_g(G)=f_n(G)$. Using this condition, we first present a short proof for
Chang and Hsieh's results~\cite{CH}. Then we investigate the $g$-extra connectivity of a particular subfamily of
HL-networks, namely, varietal hypercubes $VQ_n$. Let $n=3s+t$ with $s\geq 3$ and $0\leq t\leq 2$.
This study shows that $\k_g(VQ_n)=f_n(g)$ for $0\leq g\leq n-s$. At last, we construct a subfamily of HL-networks
to show that the inequality $\k_g(G)>f_n(g)$ can hold, and so a negative answer to Problem~A is given.

\section{Preliminaries}
Throughout this paper only undirected simple connected graphs without loops and multiple
edges are considered. Unless stated otherwise, we follow Bondy and Murty~\cite{BMBook}
for terminology and definitions.

\subsection{Fundamental graph and group terminologies}

A {\em graph} $G=(V, E)$ is comprised of a vertex set $V$ and an edge set $E$ , where $V$ is a
finite set and $E$ is a subset of $\{\{u,v\}\ |\ u,v\in V,u\neq v\}$. Two vertices, $u$ and $v$, are {\em adjacent}
if $\{u,v\}\in E$, and $u$ and $v$ are the end-vertices of $\{u,v\}$. A {\em subgraph} of $G$ is a graph $H=(V',E')$ such that $V'\subseteq V$
and $E'\subseteq E$. We use $H\subseteq G$ to denote that $H$ is a subgraph of $G$. Given a vertex set $U\subseteq V$, the subgraph of
$G$ induced by $U$ is the graph $G[U]=(U, E')$, where $E'=\{\{u,v\}\in E\ |\ u, v\in U\}$. For a set of vertices and/or edges, denoted by $F$,
in $G$, the notation $G-F$ represents a subgraph of $G$ obtained by deleting all the elements in $F$ from $G$. The components of
a graph $G$ are its maximal connected subgraphs.

A {\em $k$-path} $P_k=v_0v_1\ldots v_k$ for $k\geq 1$ in a graph $G$ is a sequence of distinct vertices such that any two consecutive
vertices are adjacent; we call $v_0$ and $v_k$ the end-vertices of the path. A {\em $k$-cycle} $C_k=(v_1, v_2, \ldots , v_k, v_1)$ for $k\geq 3$ is a sequence
of vertices in which any two consecutive vertices are adjacent, where $v_1, v_2, \ldots , v_k$ are all distinct. A complete graph $K_n$ is
a graph comprised of n pairwise adjacent vertices. A complete bipartite graph $K_{m,n}$ is a graph comprised of two partite sets of
vertices of sizes $m$ and $n$, respectively, such that two vertices are adjacent if and only if they are in different partite sets.

Let $G=(V, E)$ be a graph. The {\em neighborhood} of a vertex $u$ in a subgraph $H\subseteq G$, denoted by $N_H(u)$, is the set of all
vertices adjacent to $u$ in $H$. The degree of a vertex $u$ in $G$, denoted by $d_G(u)$, is the number of the vertices adjacent to $u$
in $G$. Note that $d_G(u)=|N_G(u)|$. For a vertex subset $V'\subseteq V$, the neighborhood of $V'$ in a subgraph $H\subseteq G$  is defined as
$N_H(V')=\bigcup_{v\in V'}(N_H(v))-V'$.

A {\em group} is a nonempty set $G$ together an binary operation $*$
defined on $G$ and satisfies the following properties:
\begin{enumerate}
  \item [{\rm (1)}]\ $*$ is associative, that means $(a*b)*c=a*(b*c)$ for any $a,b,c\in G$;
  \item [{\rm (2)}]\ $(G,*)$ has an identity $e$, that means $a*e=e*a=a$ for any $a\in G$;
  \item [{\rm (3)}]\ For any element $a\in G$, there exists an
  inverse element $a^{-1}\in G$ such that $a*a^{-1}=a^{-1}*a=e$;
\end{enumerate}

Throughout this paper, all groups are finite.
If a subset $H$ of a group $G$ is itself a group under the operation
of $G$, we say that $H$ is a {\em subgroup} of $G$ and denoted by $H\leq G$.
A subgroup $H$ of a group $G$ is called {\em normal}, denoted by $H\unlhd G$,
if $H^g=H$, $\forall g\in G$.
For a subset $S$ of a group $G$,
the intersection of all subgroups of $G$ containing $S$ is called the
{\em subgroup generated by $S$}, denoted by $\lg S\rg$.
For an element $g$ in a group $G$,
the {\em order} of $g$ is the smallest positive integer $n$ satisfying $g^n=e$.
An element of $G$ of order $2$ is also called an {\em involution}. Clearly,
if $g$ is an involution, then $g=g^{-1}$.

Let $n$ be a positive integer. Throughout this paper, $\mz_n$ represents the cyclic group
of order $n$ as well as the ring of integers modulo $n$.

An {\em isomorphism} from a simple graph $G$ to a simple graph $H$ is a bijection $\pi: V(G)\rightarrow V(H)$ such that
$\{u,v\}\in E(G)$ if and only if $\{\pi(u),\pi(v)\}\in E(H)$. If there is an isomorphism from $G$ to $H$, we
say that $G$ and $H$ are {\em isomorphic} and write $G\cong H$. An isomorphism from the graph $G$ onto itself is called
an {\em automorphism} of $G$. The set of all automorphisms of
the graph $G$, with the operation of composition, is the
automorphism group of $G$, denoted by $\Aut(G)$. We say that
$G$ is {\em vertex-transitive} if for any two vertices $u, v\in V(G)$,
there exists an automorphism $\pi\in\Aut(G)$ such that
$\pi(u)=v$.

Given a finite
group $G$ and a subset $S\subseteq G\setminus\{e\}$ such that
$S=S^{-1}=\{s^{-1}\ |\ s\in S\}$, where $e$ is the identity element of $G$,
the {\em Cayley graph} $\Cay(G,S)$ on $G$ with respect to $S$ is
defined to have vertex set $G$ and edge set $\{\{g,sg\}\mid g\in
G,s\in S\}$. A Cayley graph $\Cay(G,S)$ is connected if and only if
$S$ generates $G$. Given a $g\in G$, define the permutation $R(g)$
on $G$ by $x\mapsto xg, x\in G$. Then $R(G)=\{R(g)\ |\ g\in G\}$,
called the {\em right regular representation} of $G$, is a
permutation group isomorphic to $G$. It is well-known that
$R(G)\leq\Aut(\Cay(G,S))$. So, $\Cay(G,S)$ is vertex-transitive. In
general, a vertex-transitive graph $X$ is isomorphic to a Cayley
graph on a group $G$ if and only if its automorphism group has a
subgroup isomorphic to $G$, acting regularly on the vertex set of
$X$ (see \cite[Lemma~16.3]{B}).



\subsection{Preliminary results}

Given a positive integer $n$, let $f_n(g)=n(g+1)-\frac{g(g+3)}{2}$ be a function of $g$.\medskip

\begin{prop}{\rm \cite[Remark~2.2]{YM} or \cite[Lemma~3.1]{Zhu}}\label{fuction}
If $0\leq g\leq n-2$, then $f_n(g)$ is strictly monotonically increasing. Moreover, the
maximum of $f_n(g)$ is $f_n(n-2)=\frac{n(n-1)}{2}+1$ and $f_{n}(n-1)=f_n(n-2)>f_{n}(n)=f_n(n-3)>f_n(g)$
for $0\leq g\leq n-4$.
\end{prop}\medskip

The following result first appeared in \cite[Remark~2.2]{YM} without proof. Here we give a detailed proof.\medskip

\begin{lem} \label{inequation}
Let $0\leq g_1,g_2\leq n-2$ and $0\leq g\leq n-3$. If $g_1+1+g_2+1>g+1$, then $f_{n-1}(g_1)+f_{n-1}(g_2)\geq f_n(g)+1$.
\end{lem}\medskip

\demo We first assume that $g_1,g_2>g$. By Proposition~\ref{fuction},
$$f_{n-1}(g_1)+f_{n-1}(g_2)\geq 2f_{n-1}(g)\geq f_n(g)+f_{n-2}(0).$$
Since $n\geq 3$, one has $f_{n-2}(0)=n-2\geq1$, and so $f_{n-1}(g_1)+f_{n-1}(g_2)\geq f_n(g)+1$, as required.

Now assume that either $g_1\leq g$ or $g_2\leq g$. Without loss of generality,
let $g_1\leq g$. Then $g_2>g-g_1-1$.
Clearly, $g-g_1-1\leq g-1\leq (n-1)-3$.

If $g-g_1-1<0$, then $g_1>g-1$, and since $g_1\leq g$, one has $g_1=g$.
It follows that {\footnotesize$$f_{n-1}(g_1)+f_{n-1}(g_2)\geq f_{n-1}(g)+f_{n-1}(0)=f_n(g)+(n-g-2).$$}
As $g\leq n-3$, one has $n-g-2\geq 1$, implying $f_{n-1}(g_1)+f_{n-1}(g_2)\geq f_n(g)+1$, as required.

If $g-g_1-1\geq 0$, then by Proposition~\ref{fuction}, we have
{\footnotesize$$
\begin{array}{lll}
f_{n-1}(g_1)+f_{n-1}(g_2) &> & f_{n-1}(g_1)+f_{n-1}(g-g_1-1)\\
&=& f_n(g)+g_1(g-g_1-1)\geq f_n(g).
\end{array}
$$}
It follows that $f_{n-1}(g_1)+f_{n-1}(g_2)\geq f_n(g)+1$, as required.
\hfill\qed

The following results are very useful when used to study
the extra connectivity of HL-networks.\medskip

\begin{prop}{\rm \cite[Lemma~4]{FL}}\label{basic}
For any integer $g\geq 0$ and any integer $n\geq \lceil\frac{g+2}{2}\rceil$, for any
$X_n\in\mathbb{L}_n$ and $U\subseteq V(X_n)$ with $|U|=g+1$, we have
$|N_{X_n}(U)|\geq f_n(g)$.
\end{prop}\medskip

\begin{prop}\label{com-zhu}{\rm\cite[Lemma~4 \& Corollary~1]{Zhu2008}}
Let $X_n\in\mathbb{L}_n$. Then the girth of $X_n$ is $4$, and any two vertices of $X_n$ have at most two common neighbors.
\end{prop}\medskip

The following result is about the hypercubes.\medskip

\begin{prop}{\rm \cite[Lemma~2.1]{YM}}\label{Qnbasic}
For any integer $g\geq 0$ and any integer $n\geq 4$, for any
$U\subseteq V(Q_n)$ with $|U|=g+1$, we have
$|N_{Q_n}(U)|\geq f_n(g)$.
\end{prop}\medskip

Following Latifi~\cite{Latifi},
we express $Q_n$ as $D_0\bigodot D_1$, where $D_0$ and $D_1$ are the two
$(n-1)$-subcubes of $Q_n$ induced by the vertices with the $i$th coordinates 0 and 1,
respectively. Sometimes we use $X^{i-1}0X^{n-i}$
and $X^{i-1}1X^{n-i}$ to denote $D_0$ and $D_1$, where
$X\in\mz_2$. Clearly, the vertex $v$ in one $(n-1)$-subcube has exactly one neighbor
$v_0$ in another $(n-1)$-subcube. The following lemma presents a generalization of
\cite[Theorem~3.2]{YM}\medskip

\begin{lem}\label{star}
Let $n\geq 4$ and $g\geq 0$. Let $U\subseteq V(Q_n)$ such that $|U|=g+1$ and
$Q_n[U]$ is connected.
If $|N_{Q_n}(U)|=f_n(g)$, then $Q_n[U]$ is a star.
\end{lem}\medskip

\demo We will verify the lemma by induction on $g$. The result is clearly true for $g=0$.
We assume that $g\geq 1$ and the result is true for $h<g$. Next, we verify that this result is also true for
$h=g$. Since $|U|\geq 2$, we can take two distinct vertices, say $x=x_1x_2\ldots x_i\ldots x_n$
and $y=y_1y_2\ldots y_i\ldots y_n$ in $U$ such that $x_i=0$ and $y_i=1$. Let
$V_0=V(X^{i-1}0X^{n-i})$ and $V_1=V(X^{i-1}1X^{n-i})$. Then $x\in V_0$ and $y\in V_1$, and so
$U_i=V_i\cap U$ is non-empty for $i=0,1$. Without loss of generality, assume that $|U_0|\leq |U_1|$.
Letting $|U_0|=N$, we have $N\leq \lfloor\frac{g+1}{2}\rfloor$.
By Proposition~\ref{Qnbasic}, $|N_{V_0}(U_0)|\geq f_{n-1}(N-1)$
and $|N_{V_1}(U_1)|\geq f_{n-1}(g-N)$.
Note that $N_{V_0}(U_0)\cup N_{V_1}(U_1)\subseteq N_{Q_n}(U)$ and $N_{V_0}(U_0)\cap N_{V_1}(U_1)=\emptyset$.
It follows that $|N_{V_0}(U_0)|+|N_{V_1}(U_1)|\leq |N_{Q_n}(U)|=f_n(g)$,
and hence
\begin{equation}\label{eq-0}
\begin{array}{lll}
0&\geq &|N_{V_0}(U_0)|+|N_{V_1}(U_1)|-|N_{Q_n}(U)|\\
&\geq & f_{n-1}(N-1)+f_{n-1}(g-N)-f_n(g)\\
&=&-(N-1)(N-g).
\end{array}
\end{equation}
Remember that $1\leq N\leq \lfloor\frac{g+1}{2}\rfloor$. If $N>1$,
then $N<g$ and so $0\geq -(N-1)(N-g)>0$, a contradiction. Thus, $N=1$ and so $|U_1|=g$.
Let $U_0=U\cap V_0=\{v\}$. Since $v$ has only one neighbor in $V_1$,
$v$ is a pendant vertex of $Q_n[U]$, and so $Q_n[U_1]$ is connected.
Since $N=1$, by the above equation (\ref{eq-0}),
we have $|N_{V_1}(U_1)|=f_{n-1}(g-1)$ and $|N_{Q_n}(U)|=|N_{V_0}(U_0)|+|N_{V_1}(U_1)|$.
It follows that $|N_{Q_n}(U_1)|=|N_{V_1}(U_1)|+|N_{V_0}(U_1)|=f_{n-1}(g-1)+g=f_n(g-1)$.
By the induction hypothesis,
we have $Q_n[U_1]$ is a star. Assume that $U_1=\{u,u_1,\ldots,u_{g-1}\}$
and $E(Q_n[U_1])=\{\{u,u_i\}\ |\ 1\leq i\leq g-1\}$.
If $v$ is adjacent to $u$, then $Q_n[U]\cong K_{1,g}$.
Assume that $v$ is adjacent to some $u_i$.
If $g\leq 2$, then $Q_n[U]$ must be a star. Suppose $g>2$.
If there is a $w\in N_{V_0}(U_1-u_i)\setminus N_{V_0}(v)$,
then $|N_{V_0}(U_0)|+|N_{V_1}(U_1)|=|N_{Q_n}(U)|\geq 1+|N_{V_0}(U_0)|+|N_{V_1}(U_1)|$, a contradiction.
Thus, $N_{V_0}(U_1-u_i)\subseteq N_{V_0}(v)$. Take $u_j\in U_1-\{u,u_i\}$. Then
the neighbor $w_j$ of $u_j$ in $V_0$ is adjacent to $v$.
This implies that $(v,u_i,u,u_j,w_j)$ is a cycle of length $5$,
contrary to the fact that $Q_n$ is bipartite.
\hfill\qed

\section{Lower bound on $g$-extra connectivity of HL-networks}


A graph is said to be {\em hyper-$\k_g$} if the deletion of each minimum $R_g$-cutset
creates exactly two components, one of which has $g+1$ vertices.
Clearly, a hyper-$\k_0$ graph is also hyper-$\k$ (see \cite{Meng2003} for the definition of hyper-$\k$).
The following lemma shows that every HL-network is hyper-$\k$.

\medskip

\begin{lem}\label{super-k}
For any $X_n\in {\mathbb L}_n$, $X_n$ is hyper-$\k$ for $n\geq 2$.
\end{lem}\medskip

\demo We will prove the lemma by induction on $n$.
Let $S$ be a minimum vertex-cut of $X_n$.
By \cite{FH}, we have $|S|=n$.
The result is clearly true when $n=2$. In what follows,
assume that $n\geq 3$, and that the result holds for $X_{n-1}$.

Suppose $X_n=X_{n-1}^0\oplus X_{n-1}^1$. Let $S_i=S\cap V(X_{n-1}^i)$
with $i=0,1$. Since $n\geq 3$, one has $2^{n-1}-n\geq 1$, and so
there is at least one edge
between $X_{n-1}^0-S_0$ and $X_{n-1}^1-S_1$.
Without loss of generality, assume that $|S_0|\leq |S_1|$.
If $|S_0|=0$, then $X_{n-1}^0$ is connected, and since
each vertex of $X_{n-1}^1$ has a neighbor in $X_{n-1}^0$,
it follows that $X_n-S$ is connected, a contradiction. Thus, $|S_0|\geq 1$.
If $|S_0|>1$, then $|S_0|\leq |S_1|<n-1$, and so $X_{n-1}^i-S_i$ is connected
for $i=0,1$. Since there is at least one edge
between $X_{n-1}^0-S_0$ and $X_{n-1}^1-S_1$,
$X_n-S$ is still connected, a contradiction. Consequently, $|S_0|=1$,
$|S_1|=n-1$ and $X_{n-1}^1-S_1$ is disconnected. By the induction hypothesis,
$X_{n-1}^1-S_1$ has exactly two components, one of which is a singleton, say $u$.
Clearly, $X_{n-1}^0-S_0$ is connected. If $n\geq 4$, then $2^{n-1}-n\geq 4$, and so there are at
least four edges between $X_{n-1}^0-S_0$ and $X_{n-1}^1-S_1$, the component of $X_{n-1}^1-S_1$
which is not the singleton is connected to $X_{n-1}^0-S_0$.
It follows that $X_n-S$ has exactly two components, one of which is the singleton $u$, as required.
If $n=3$, then $X_{n-1}^1-S_1$ is a null graph with two vertices, say $u,v$. Since there is at least one edge
between $X_{n-1}^0-S_0$ and $X_{n-1}^1-S_1$, either $u$ or $v$ is connected to $X_{n-1}^0-S_0$.
Again, $X_n-S$ has exactly two components, one of which is a singleton, as required. \hfill\qed

As a slight generalization of \cite[Theorem~3.3]{Zhu}, the following result
shows that an $n$-dimensional HL-network with at most $f_n(g)$ faulty vertices, where $0\leq g\leq n-3$,
has a very large component. \medskip

\begin{lem}\label{structure}
Let $n\geq 5$, $0\leq g\leq n-3$. For any $X_n\in {\mathbb L}_n$,
let $S\subseteq V(X_n)$. If $|S|\leq f_n(g)-k$ with $0\leq k\leq 1$,
then $X_n-S$ has a component with at least $2^n-|S|-(g+1-k)$ vertices.
\end{lem}\medskip

\demo We will prove the lemma by induction on $g$. Let $g=0$. Then
$|S|\leq n-k$. By Lemma~\ref{super-k}, if $k=1$, then $G-S$ is connected,
and if $k=0$ and $|S|=n$, then $G-S$ either is connected, or has two components, one of which is a singleton.
This implies that the result is true for $g=0$.

In what follows,
assume that $g\geq 1$, and that the result holds for $g-1$.
We shall verify that it is also true for $g$.
Suppose $X_n=X_{n-1}^0 \oplus X_{n-1}^1$. Let $S_i=S\cap V(X_{n-1}^i)$
with $i=0,1$. Without loss of generality, assume that $|S_0|\leq |S_1|$.
\medskip

\f{\bf Case~1}\ $|S_0|\leq n-2$

In this case, $X_{n-1}^0-S_0$ is connected because $\kappa(X_{n-1}^0)=n-1$ by Lemma~\ref{super-k}.
Let $C$ be the component of $X_n-S$ containing $X_{n-1}^0-S_0$, and let
$C'$ be the union of all other components of $X_n-S$. Then $C'\subseteq X_{n-1}^1$.
Noting that the edges between $V(X_{n-1}^0)$ and $V(X_{n-1}^1)$ form a perfect matching, one has
$|V(C')|=|N_{X_{n-1}^0}(V(C'))|\leq |S_0|\leq n-2$. Clearly,
$|N_{X_n}(V(C'))|\leq |S|$. By Proposition~\ref{basic}, we have $|N_{X_n}(V(C'))|\geq f_n(|V(C')|-1)$.
It follows that $f_n(g)-k\geq|S|\geq |N_{X_n}(V(C'))|\geq f_n(|V(C')|-1)$.
By Proposition~\ref{fuction}, if $k=0$, then $|V(C')|-1\leq g$ and so $|V(C')|\leq g+1$;
if $k=1$, then $|V(C')|-1<g$ and so $|V(C')|\leq g=(g+1)-1$. Thus, we always have
$|V(C')|\leq (g+1)-k$ and so $|V(C)|\geq 2^n-|S|-(g+1-k)$. \medskip


\f{\bf Case~2}\ $|S_0|>n-2$

In this case, we have $$
\begin{array}{lll}
|S_0|\leq |S_1|&=&|S|-|S_0|\\
&\leq & f_n(g)-k-(n-1)\\
&=&f_{n-1}(g-1)-k.
\end{array}
$$
By the induction hypothesis, for each $i=0,1$,
$X_{n-1}^i$ has a component $A_i$ with at least $2^{n-1}-|S_i|-(g-k)$ vertices.
Clearly, $|S|+2g\leq f_n(g)-k+2g$.
By Proposition~\ref{fuction}, $f_n(g)\leq f_n(n-3)=\frac{n^2-n}{2}$. It follows that
$|S|+2g\leq f_n(n-3)-k+2(n-3)=\frac{n^2+3n-12-2k}{2}<2^{n-1}$ because $n\geq 5$.
In view of the fact that there are $2^{n-1}$ vertex disjoint edges between
$V(X_{n-1}^0)$ and $V(X_{n-1}^1)$, there exist some edges between $A_0$ and
$A_1$ in $X_n-F$. Let $C$ be the component of $X_n-S$ containing
$A_1$ and $A_2$, and let $C'$ be the union of all other components of $X_n-S$.
Let $V_i=V(C')\cap V(X_{n-1}^i)$ for $i=0,1$. Then $|V_0|,|V_1|\leq g-k\leq n-3$.
Note that $|S|\geq |N_{X_{n-1}^0}(V_0)|+|N_{X_{n-1}^1}(V_1)|$.
By Proposition~\ref{basic}, $|N_{X_{n-1}^0}(V_0)|+|N_{X_{n-1}^1}(V_1)|\geq f_{n-1}(|V_0|-1)+f_{n-1}(|V_1|-1)$.
It follows that $|S|\geq f_{n-1}(|V_0|-1)+f_{n-1}(|V_1|-1)$.
Since $|S|\leq f_n(g)-k$, one has $f_n(g)-k\geq f_{n-1}(|V_0|-1)+f_{n-1}(|V_1|-1)$.
From Lemma~\ref{inequation} it follows that $|V_0|+|V_1|\leq g+1$.
If $k=0$, then $|V_0|+|V_1|\leq g+1=g+1-k$. Let $k=1$. Suppose that
$|V_0|+|V_1|=g+1$. Letting $g_0=|V_0|-1$, we have
{\footnotesize$$
\begin{array}{l}
f_{n-1}(|V_0|-1)+f_{n-1}(|V_1|-1)=f_n(g)+g_0[g-g_0-1].
\end{array}
$$}
As $g-g_0-1=|V_1|-1\geq 0$, one has $f_{n-1}(|V_0|-1)+f_{n-1}(|V_1|-1)\geq f_n(g)$,
contrary to the fact that $f_n(g)-1\geq f_{n-1}(|V_0|-1)+f_{n-1}(|V_1|-1)$.
Thus, $|V_0|+|V_1|\leq g$.
As a result, we always have $V(C')=|V_0|+|V_1|\leq g+1-k$ and so $|V(C)|\geq 2^n-|S|-(g+1-k)$, as required. \hfill\qed

By Lemma~\ref{structure}, we immediately have the following result.\medskip

\begin{theorem}\label{extra-0}
Let $n\geq 5$, $0\leq g\leq n-3$. For any $X_n\in {\mathbb L}_n$,
let $S$ be a $R_g$-cutset of $X_n$. If $|S|\leq f_n(g)=n(g+1)-\frac{g(g+3)}{2}$,
then $|S|=f_n(g)$ and $X_n$ is hyper-$\k_g$.
In particular, $\k_g(X_n)\geq f_n(g)$.
\end{theorem}
\medskip

\demo If $|S|\leq f_n(g)-1$ then by Lemma~\ref{structure}, $X_n-S$ would have a component
with at most $g$ vertices, contrary to the fact that $S$ is a $R_g$-cutset.
Therefore, $|S|=f_n(g)$. Again, by Lemma~\ref{structure}, $X_n-S$ has a component, say $C$, with at least
$2^n-|S|-(g+1)$ vertices. Let $A$ be a component of $X_n-S$ different from $S$.
Then $|V(A)|\leq g+1$. Since $S$ is a $R_g$-cutset, we must have $|V(A)|=g+1$. Therefore,
$X_n-S$ has exactly two components that are $C$ and $A$. \hfill\qed

\f{\bf Remark~1}\ After this work was finished, it came to our notice that Yang and Liu~\cite[Remark~4.11]{YL}
also proved that for $0\leq g\leq n-4$, if $X_n$ has a $R_g$-cutset $S$ with $|S|\leq f_n(g)$,
then $\k_g(X_n)=f_n(g)$.
With our approach, however, we were able to further obtain that $X_n$ is also hyper-$\k_g$.

\section{Some HL-networks $X_n$ with $\k_g(X_n)=f_n(g)$}
In this section, we shall show that the lower bound on $\k_g(X_n)$ in Theorem~\ref{extra-0}
is best possible. We first give a sufficient condition for the equality $\k_g(X_n)=f_n(g)$.\medskip

\begin{theorem}\label{extra-1}
Let $n\geq 5$, $0\leq g\leq n-3$, and $X_n\in {\mathbb L}_n$.
If $X_n$ has a connected subgraph, say $A$, such that $|V(A)|=g+1$ and
$|N_{X_n}(A)|=f_n(g)$, then $\k_g(X_n)=f_n(g)$.
\end{theorem}\medskip

\demo By Theorem~\ref{extra-0}, it suffices to show that $\k_g(X_n)\leq f_n(g)$.
Set $S=N_{X_n}(A)$. Then $A$ is a component of $G-S$ with
$g+1$ vertices. Since $|S|=f_n(g)$, from Lemma~\ref{structure} it follows
that $G-S$ has a component, say $A'$, with at least $2^n-|S|-(g+1)$ vertices.
By Proposition~\ref{fuction}, $2^n-f_n(g)-(g+1)\geq 2^n-f_n(n-3)-(n-2)$.
Since $n\geq 5$, one has $2^n-f_n(n-3)-(n-2)=2^n-\frac{n^2+n-4}{2}>n\geq g+3$.
It follows that $|V(A')|\geq 2^n-|S|-(g+1)=2^n-f_n(g)-(g+1)>g+3$, and so
$A'\neq A$. Then $V(X_n)=V(A)\cup V(A')\cup S$, which implies that
$S$ is a $g$-extra vertex cut of $X_n$. Hence, $\k_g(X_n)\leq f_n(g)$.\hfill\qed


This theorem is very powerful when used to determine the $g$-extra connectivity of
HL-networks for small $g$.
For example, we can use it to give a short proof of \cite[Theorems~1,2]{CH}.\medskip

\begin{cor}
For any $X_n\in\mathbb{L}_n$,
$\k_2(X_n)=3n-5$ for $n\geq 5$ and $\k_3(X_n)=4n-9$ for $n\geq 6$.
\end{cor}\medskip

\demo From the definition of HL-networks,  $X_n$ has a subgraph isomorphic to a
$k$-dimensional $HL$-network for each $1\leq k\leq n$. Pick a subgraph, say $C$, isomorphic to a $2$-dimensional
HL-network. Then $C$ is a $4$-cycle. Let $C=(u_0,u_1,u_2,u_3,u_0)$ and let $A=u_0u_1u_2$.
By Proposition~\ref{com-zhu}, any two vertices of $X_n$ have at most two common neighbors.
Then $|N_{X_n}(A)|=3n-5=f_n(2)$. Since $n\geq 5$, one has $g=2\leq n-3$.
From Theorem~\ref{extra-1} it follows that $\k_2(X_n)=3n-5$.

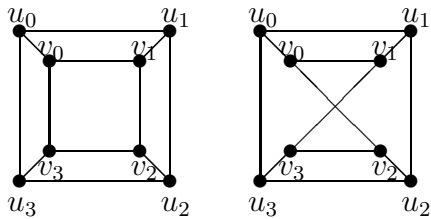
\begin{figure}[ht]
\begin{center}
\unitlength 4mm
\begin{picture}(15,6)

\put(1, 1){\circle*{0.4}} \put(0.6, 0){$u_3$}

\put(1,1){\line(1, 0){5}} \put(1,1){\line(0, 1){5}}

\put(6, 1){\circle*{0.4}}\put(5.7, 0){$u_2$}

\put(6,1){\line(0, 1){5}}\put(1,6){\line(1, 0){5}}

\put(1, 6){\circle*{0.4}} \put(0.6, 6.3){$u_0$}

\put(6, 6){\circle*{0.4}}\put(5.7, 6.3){$u_1$}


\put(2, 2){\circle*{0.4}} \put(1.6, 1.2){$v_3$}

\put(2,2){\line(1, 0){3}} \put(2,2){\line(0, 1){3}}

\put(2,2){\line(-1, -1){1}}\put(5,2){\line(1, -1){1}}

\put(5, 2){\circle*{0.4}}\put(4.7, 1.2){$v_2$}

\put(5,2){\line(0, 1){3}}\put(2,5){\line(1, 0){3}}

\put(2, 5){\circle*{0.4}} \put(1.6, 5.2){$v_0$}

\put(2,5){\line(-1, 1){1}}

\put(5, 5){\circle*{0.4}}\put(4.7, 5.2){$v_1$}

\put(5,5){\line(1, 1){1}}


\put(9, 1){\circle*{0.4}} \put(8.6, 0){$u_3$}

\put(9,1){\line(1, 0){5}} \put(9,1){\line(0, 1){5}}

\put(14, 1){\circle*{0.4}}\put(13.7, 0){$u_2$}

\put(14,1){\line(0, 1){5}}\put(9,6){\line(1, 0){5}}

\put(9, 6){\circle*{0.4}} \put(8.6, 6.3){$u_0$}

\put(14, 6){\circle*{0.4}}\put(13.7, 6.3){$u_1$}


\put(10, 2){\circle*{0.4}} \put(9.6, 1.2){$v_3$}

\put(10,2){\line(1, 0){3}} \put(10,2){\line(1, 1){3}}

\put(10,2){\line(-1, -1){1}}\put(13,2){\line(1, -1){1}}

\put(13, 2){\circle*{0.4}}\put(12.7, 1.2){$v_2$}

\put(13,2){\line(-1, 1){3}}\put(10,5){\line(1, 0){3}}

\put(10, 5){\circle*{0.4}} \put(9.6, 5.2){$v_0$}

\put(10,5){\line(-1, 1){1}}

\put(13, 5){\circle*{0.4}}\put(12.7, 5.2){$v_1$}

\put(13,5){\line(1, 1){1}}
\end{picture}
\end{center}\vspace{-.5cm}
\caption{Two $3$-dimensional HL-networks} \label{Fig-3.1}
\end{figure}

Now take a subgraph, say $H$, isomorphic to a $3$-dimensional HL-network.
Then $H$ is one of the following two graphs in Figure~\ref{Fig-3.1}.
If $H$ is the first graph, then let $A$ be the star $K_{1,3}$ with vertices $u_0,v_0,u_1,u_3$.
If $H$ is the second graph, then let $A$ be the $3$-path $u_0v_0v_1v_3$.
For both cases, it is also easy to check that $|N_{X_n}(A)|=4n-9=f_n(3)$.
Since $n\geq 6$, one has $g=3\leq n-3$.
By Theorem~\ref{extra-1}, we have $\k_3(X_n)=4n-9$.\hfill\qed

In what follows, we shall consider the $g$-extra connectivity of the
varietal hypercubes which were proposed by Cheng and Chuang~\cite{CC}.
An $n$-dimensional varietal hypercube, denoted by $VQ_n$, is defined
recursively as follows.\medskip

\begin{defi}
The $VQ_1$ is the complete graph of two vertices labeled with $0$ and $1$,
respectively. Assume that $VQ_{n-1}$ has been constructed for $n\geq 1$.
Let $VQ_n^0 ($resp. $VQ_n^1)$ be the graph obtained from $VQ_{n-1}$ be inserting a
$0 ($resp. $1)$ in front of each vertex-labeling in $VQ_{n-1}$.
Then the $VQ_n$ is obtained by joining vertices in $VQ_n^0$ and $VQ_n^1$,
according to the rule: a vertex
$0x_{n-1}x_{n-2}x_{n-3}\ldots x_2x_1$ in $VQ_n^0$ and a vertex
$0y_{n-1}y_{n-2}y_{n-3}\ldots y_2y_1$ in  $VQ_n^1$ are adjacent in
$VQ_n$ if and only if one of the following holds:
\begin{enumerate}
  \item [{\rm (1)}]\ $x_{n-1}x_{n-2}x_{n-3}\ldots x_2x_1=y_{n-1}y_{n-2}y_{n-3}\ldots y_2y_1$ if $3\ |\ n$;
  \item [{\rm (2)}]\ $x_{n-3}\ldots x_2x_1=y_{n-3}\ldots y_2y_1$ and $(x_{n-1}x_{n-2},$ $y_{n-1}y_{n-2})\in$ $\{(00,00),$ $(01,01),$ $ (10,11),$ $(11,10)\}$ if $3\ |\ n$.
\end{enumerate}
 \end{defi}\medskip

Clearly, for $n\geq 1$, $VQ_n$ is an $n$-regular graph with vertex set
$V=\{x_nx_{n-1}\ldots x_2x_1\ |\ x_i\in\mz_2, 1\leq i\leq n\}$.
In what follows, we always assume that $n=3s+t$, where $0\leq t\leq 2$ and $s$ is a non-negative integer.\medskip

\begin{lem}\label{adjacency}
Let $u=x_nx_{n-1}\ldots x_3x_2x_1\in V(VQ_n)$. Then a vertex $v=y_{n}y_{n-1}\ldots y_3y_2y_1$
is adjacent to $u$ if and only if one of the following holds:
\begin{enumerate}
  \item [{\rm (1)}]\ for some $1\leq i\leq n$ with $3\nmid i$, $y_i=\overline{x_i}$, and for all $1\leq j\leq n$ with $j\neq i$,
  $y_j=x_j$;
  \item [{\rm (2)}]\ for some $1\leq i\leq n$ with $3\ |\ i$,  $y_i=\overline{x_i}, y_{i-2}=x_{i-1}+x_{i-2}$,
   and for all $1\leq j\leq n$ with $j\neq i,i-2$,  $y_j=x_j$.
\end{enumerate}
\end{lem}\medskip

\demo We shall first prove the sufficiency by using induction on $n$.
The result is obviously true for $n=1$ because $VQ_1$ is a complete graph with two vertices $\{0,1\}$.
Now assume that the result holds for $n-1$.
If $i<n$, then by induction, $x_{n-1}\ldots x_3x_2x_1$ is adjacent to $y_{n-1}\ldots y_3y_2y_1$ in $VQ_{n-1}$.
Since $x_n=y_n$, either $u,v\in V(VQ_{n}^0)$ or $u,v\in V(VQ_{n}^1)$. It follows that $u$ and $v$ are also adjacent in $VQ_n$.
If $i=n$, we may assume that $x_n=0$ and $y_n=1$. If $3\nmid n$, then $y_j=x_j$ for all $0\leq j\leq n-1$,
and then by the definition of $VQ_n$, $u$ and $v$ are adjacent.
If $3\ |\ n$, then $y_{n-2}=x_{n-1}+x_{n-2}$, and $y_j=x_j$ for all $1\leq j\leq n$ with $j\neq n,n-2$.
Noting that $x_i\in\mz_2$ for all $1\leq i\leq n$,
it is easy to check that $(x_{n-1}x_{n-2},y_{n-1}y_{n-2})\in\{(00,00), (01,01), (10,11), (11,10)\}$.
By the definition of $VQ_n$, $u$ and $v$ are adjacent.

For the necessity, let $N$ be the set of vertices whose coordinates satisfy the condition $(1)$ or $(2)$. 
By the sufficiency, each vertex in $N$ is adjacent to $u$. Clearly, $|N|=n$, so $N$ is just the neighborhood of $u$ in $VQ_n$
because $VQ_n$ has valency $n$. So, $v\in N$. \hfill\qed 

Remember that $n=3s+t$ with $s\geq 0$ and $0\leq t\leq 2$. For $1\leq i\leq s$, let $H_i$ be
the dihedral group of order $8$ defined as follows:
$$
\begin{array}{lll}
H_i&=&\lg a_i,b_i\ |\ a_i^4=b_i^2=1, b_i^{-1}a_ib_i=a_i^{-1}\rg.
\end{array}
$$
Let $\lg c_1\rg\times\ldots\times\lg c_t\rg\cong\underbrace{\mz_2\times\ldots\times\mz_2}_{t\ {\rm times}}$.
Set $G_n=H_1\times H_2\times\ldots\times H_s\times \lg c_1\rg\times\ldots\times\lg c_t\rg$,
and $\Omega_n=\{a_i^2,b_i,a_ib_i,c_1, \ldots, c_t\ |\ 1\leq i\leq s\}$.
Let $\Delta_n=\Cay(G_n,\Omega_n)$. 

\medskip
\begin{theorem}\label{isomorphism}
$\Delta_n\cong VQ_n$.
\end{theorem}\medskip

\demo We define a map $f$ from $V(\Delta_n)$ to $V(VQ_n)$ as follows:
{\footnotesize$$
\begin{array}{l}
(a_1^2)^{x_1}b_1^{x_2}(a_1b_1)^{x_3}\ldots  (a_s^2)^{x_{3s-2}}b_s^{x_{3s-1}}(a_sb_s)^{x_{3s}}c_1^{x_{3s+1}}\ldots c_t^{x_{3s+t}}\\
\hspace{5mm}\mapsto x_{3s+t}\ldots x_{3s+1}x_{3s}x_{3s-1}x_{3s-2}\ldots  x_3x_2x_1,
\end{array}
$$}
where $x_i\in\mz_2$ with $1\leq i\leq n$.
It is easy to check that $f$ is a bijection. Take an edge, say $\{u,v\}$, of
$\Delta_n$. Without loss of generality, let {\footnotesize$$u=(a_1^2)^{x_1}b_1^{x_2}(a_1b_1)^{x_3}\ldots  (a_s^2)^{x_{3s-2}}b_s^{x_{3s-1}}(a_sb_s)^{x_{3s}}c_1^{x_{3s+1}}\ldots c_t^{x_{3s+t}}.$$} Then $v=gu$ for some $g\in\Omega_n$.

If $g=c_i$ for some $0\leq i\leq t$, then {\footnotesize$$
\begin{array}{l}v=(a_1^2)^{x_1}b_1^{x_2}(a_1b_1)^{x_3}\ldots  (a_s^2)^{x_{3s-2}}b_s^{x_{3s-1}}(a_sb_s)^{x_{3s}}c_1^{x_{3s+1}}\ldots
c_i^{x_{3s+i}+1}\ldots  c_t^{x_{3s+t}},\end{array}$$}
and so {\footnotesize$v^f=x_{3s+t}\ldots\overline{ x_{3s+i}}\ldots x_{3s+1}x_{3s}x_{3s-1}x_{3s-2}\ldots  x_3x_2x_1$}.
By Lemma~\ref{adjacency}, $\{u^f,v^f\}\in E(VQ_n)$.

If $g=a_i^2$ for some $1\leq i\leq s$, then
{\footnotesize$$
\begin{array}{l}v=(a_1^2)^{x_1}\ldots  (a_i^2)^{x_{3i-2}+1}(b_i)^{x_{3i-1}}(a_ib_i)^{x_{3i}}\ldots
b_s^{x_{3s-1}}(a_sb_s)^{x_{3s}}c_1^{x_{3s+1}}\ldots c_t^{x_{3s+t}},\end{array}$$}
and so {\footnotesize$v^f=x_{3s+t}\ldots x_{3s+1}x_{3s}\ldots {x_{3i}}{x_{3i-1}}\overline{x_{3i-2}}\ldots  x_2x_1$}.
By Lemma~\ref{adjacency}, $\{u^f,v^f\}\in E(VQ_n)$.

If $g=b_i$ for some $1\leq i\leq s$, then {\footnotesize$$\begin{array}{l}
v=(a_1^2)^{x_1}\ldots  (a_i^2)^{x_{3i-2}}(b_i)^{x_{3i-1}+1}(a_ib_i)^{x_{3i}}\ldots
b_s^{x_{3s-1}}(a_sb_s)^{x_{3s}}c_1^{x_{3s+1}}\ldots c_t^{x_{3s+t}},\end{array}$$}
and so {\footnotesize$v^f=x_{3s+t}\ldots x_{3s+1}x_{3s}\ldots x_{3i}\overline{x_{3i-1}}x_{3i-2}\ldots  x_2x_1$}.
Again, by Lemma~\ref{adjacency}, $\{u^f,v^f\}\in E(VQ_n)$.

If $g=a_ib_i$ for some $1\leq i\leq s$, then
{\footnotesize$$\begin{array}{lll}
v&=&g(a_1^2)^{x_1}\ldots  (a_i^2)^{x_{3i-2}}(b_i)^{x_{3i-1}}(a_ib_i)^{x_{3i}}\ldots b_s^{x_{3s-1}}(a_sb_s)^{x_{3s}}c_1^{x_{3s+1}}\ldots c_t^{x_{3s+t}}\\
&=& (a_1^2)^{x_1}\ldots  (a_i^2)^{x_{3i-2}}g(b_i)^{x_{3i-1}}(a_ib_i)^{x_{3i}}\ldots b_s^{x_{3s-1}}(a_sb_s)^{x_{3s}}c_1^{x_{3s+1}}\ldots c_t^{x_{3s+t}}\\
&=& (a_1^2)^{x_1}\ldots  (a_i^2)^{x_{3i-2}+x_{3i-1}}(b_i)^{x_{3i-1}}(a_ib_i)^{x_{3i}+1}\ldots b_s^{x_{3s-1}}(a_sb_s)^{x_{3s}}c_1^{x_{3s+1}}\ldots c_t^{x_{3s+t}}.
\end{array}
$$}
and so {\footnotesize$v^f=x_{3s+t}\ldots x_{3s+1}x_{3s}\ldots\overline{x_{3i}}x_{3i-1}$ $(x_{3i-1}+x_{3i-2})$ $\ldots  x_2x_1$.}
Again, by Lemma~\ref{adjacency}, $\{u^f,v^f\}\in E(VQ_n)$.

Now we see that $f$ is an isomorphism from $\Delta_n$ to $VQ_n$. Therefore, $VQ_n\cong \Delta_n$.\hfill\qed

From Theorem~\ref{isomorphism}, we see that $VQ_n$ is a Cayley graph, and so it is vertex-transitive.\medskip

\begin{cor}{\rm \cite[Theorem~2.5]{XCX}}
$VQ_n$ is vertex-transitive.
\end{cor}\medskip

\begin{theorem}
Let $n=3s+t$ with $s\geq 3$ and $0\leq t\leq 2$. If $0\leq g\leq n-s$,
then $\k_g(VQ_n)=f_n(g)$.
\end{theorem}\medskip

\demo By Theorem~\ref{isomorphism}, $VQ_n\cong\Delta_n$, and so $\k_g(VQ_n)=\k_g(\Delta_n)$.
Set $\Omega'=\{a_i^2,$ $b_i,c_1,$ $\ldots, c_t$ $ \ |\ 1\leq i\leq s\}.$ Clearly, $|\Omega'|=n-s$.
Since $0\leq g\leq n-s$, we can take a subset $V'$ of $V(\Delta_n)$ such that $|V'|=g+1$, $e\in V'$ and $V'-\{e\}\subset \Omega'$,
where $e$ is the identity element of the group $G_n$. Clearly, $\Delta_n[V']\cong K_{1,g}$.
Note that the elements in $\Omega'$
are pair-wise commutative. So, for any distinct $a,a'\in \Omega'$, $e$ and $aa'$
are two distinct common neighbors of $a$ and $a'$. From Proposition~\ref{com-zhu} it follows that
any two elements in $\Omega'$
have exactly two common neighbors in $\Delta_n$. With an easy calculation,
we see that $|N_{\Delta_n}(V')|=f_n(g)$.
By Theorem~\ref{extra-1}, we have $\k_g(\Delta_n)=f_n(g)$. \hfill\qed

\section{A negative answer to Problem~A}

For an $n$-dimensional HL-network $X_n$, if $0\leq g\leq n-3$ and $n\geq 5$, then we
have $\k_g(X_n)\geq f_n(g)$ by Theorem~\ref{extra-0}.
In this section, we shall construct a class of HL-networks with
$g$-extra connectivity greater than $f_n(g)$.

Let $k,\ell$ be two non-negative integers.
Let $\lg c_1\rg\times\ldots\times\lg c_\ell\rg\cong\underbrace{\mz_2\times\ldots\times\mz_2}_{\ell\ {\rm times}}$,
and for $1\leq i\leq k$, let $H_i$ be
the dihedral group of order $8$ defined as follows:
$$
\begin{array}{lll}
H_i&=&\lg a_i,b_i\ |\ a_i^4=b_i^2=1, b_i^{-1}a_ib_i=a_i^{-1}\rg.
\end{array}
$$
Set $G_{k,\ell}=H_1\times H_2\times\ldots\times H_k\times \lg c_1\rg\times\ldots\times\lg c_\ell\rg$.
Let $\G_{k,\ell}$ be the Cayley graph $\Cay(G_{k,\ell},\Omega_{k,\ell})$,
where $\Omega_{k,\ell}=\{a_i^2,b_i,a_ib_i,c_1, \ldots, c_\ell\ |\ 1\leq i\leq k\}$.
Let $\mathbb{G}_n=\{G_{k,\ell}\ |\ 3k+\ell=n\}$ and let $\mathcal{G}_n=\{\G_{k,\ell}\ |\ 3k+\ell=n\}$.

\medskip
\begin{lem}
Each graph in $\mathcal{G}_{n}$ is isomorphic to an $n$-dimensional HL-network.
\end{lem}\medskip

\demo We shall verify the result by using induction on $n$.
If $n=1$, then $k=0$ and $\ell=1$, and so $\mathcal{G}_1=\{\G_{0,1}\}$.
Clearly, $\G_{0,1}=\Cay(\lg c_1\rg,\{c_1\})\cong K_2$, which belongs to $\mathbb{L}_1$.
In what follows, assume that $n>1$ and that the result is true for $h\leq n-1$.
We will verify that it is also true for $h=n$.

Let $k=0$. Then $\mathcal{G}_n=\mathcal{G}_\ell=\{\G_{0,\ell}\}$.
Note that $G_{0,\ell}\cong\underbrace{\mz_2\times\ldots\times\mz_2}_{\ell\ {\rm times}}$,
and $\Omega_{0,\ell}=\{c_1, \ldots, c_\ell\}$.
It is easy to check that $\G_{0,\ell}=\Cay(G_{0,\ell},\Omega_{0,\ell})\cong Q_\ell$
which belongs to $\mathbb{L}_\ell$.

Let $k>0$. Take any $\G_{k,\ell}\in \mathbb{G}_n$. Then $\G_{k,\ell}=\Cay(G_{k,\ell},\Omega_{k,\ell})$.
Set $M=H_1\times \ldots\times H_{k-1}\times \lg a_k^2\rg\times \lg a_kb_k\rg\times
\lg c_1\rg\times\ldots\times\lg c_\ell\rg$. Clearly, $M$ is a subgroup of $G_{k,\ell}$ of index $2$, and
$M\cong G_{k-1,\ell+2}$. Set $A=\Cay(M,\Omega_{k,\ell}-\{b_k\})$. It is easy to see that $A\cong \G_{k-1,\ell+2}\in\mathcal{G}_{n-1}$.
By the induction hypothesis, $A$ is isomorphic to an $(n-1)$-dimensional HL-network.
Note that $M$ has two cosets in $G_{k,\ell}$ that are $M$ and $Mb_k$. In view of the fact that
the map $R(b_k): v\mapsto vb_k, \forall v\in G_{k,\ell}$ is an automorphism of $\G_{k,\ell}$, it follows that
$\G_{k,\ell}[Mb_k]\cong A$, and hence it is also
isomorphic to an $(n-1)$-dimensional HL-network.
By the definition of Cayley graph, we see that for any $u\in M$,
$b_ku$ is the unique neighbor of $u$ in $Mb_k$.
This implies that the edges between $M$ and $Mb_k$ form a perfect matching.
Therefore, $\G_{k,\ell}$ is isomorphic to an $n$-dimensional HL-network.\hfill\qed

\begin{lem}\label{common neighbor}
For any $\G_{k,\ell}\in \mathcal{G}_n$, and for any $u,w\in N_{\G_{k,\ell}}(v)$,
if $u$ and $w$ have a unique common neighbor, say $v$, then
$\{u,w\}=\{b_iv,a_ib_iv\}$ for some $1\leq i\leq k$.
\end{lem}\medskip

\demo Since $v$ is a common neighbor of $u$ and $w$, by the definition of Cayley graph, we have
$u=gv$ and $w=g'v$ for some $g,g'\in\Omega_{k,\ell}$. If $g$ commutes with
$g'$, then $gg'v=g'gv$ is also a common neighbor of $u$ and $w$. The uniqueness of $v$
implies that $gg'v=v$, and so $gg'=e$, where $e$ is the identity element of the group $G_{k,\ell}$.
Note that all elements of $\Omega_{k,\ell}$ are involutions (elements of order $2$). It follows that $g=g'$, forcing that
$u=gv=g'v=w$, a contradiction.
Thus, $g$ does not commute with $g'$. By the structure of the group $G_{k,\ell}$,
we see that the only possibility is $\{g,g'\}=\{b_i,a_ib_i\}$ for some $1\leq i\leq k$,
and so $\{u,w\}=\{b_iv,a_ib_iv\}$.\hfill\qed

\begin{lem}\label{component}
Let $n\geq 5$ and let $0\leq g\leq n-4$. For any $\G_{k,\ell}\in \mathcal{G}_n$, let $A$ be a connected subgraph of
$\G_{k,\ell}$ such that $|V(A)|=g+1$.
If $|N_{\G_{k,\ell}}(V(A))|=f_n(g)$, then $A$ is isomorphic either to $K_{1,g}$ or to $P_3$.
\end{lem}\medskip

\demo 
We shall prove the theorem by using induction on $g$.
Obviously, the result is true for $g=0$. In what follows, assume that
$g>0$ and the result is true for $h<g$. We shall verify that it
is also true for $h=g$.

By the vertex-transitivity of $\G_{k,\ell}$,
we may assume that the identity element $e$ of the group
$G_{k,\ell}$ is contained in $V(A)$. We consider the following two cases:\medskip

\f{\bf Case~1}\ For any edge $\{x,y\}$ of $A$, $xy^{-1}\in\{a_i^2\ |\ 1\leq i\leq k\}$.

For any vertex, say $v$, of $A$, $A$ has a path with end vertices $e$ and $v$ because of the connectivity of $A$.
Assume that this path is $v_0v_1v_2\ldots v_j$ with $v_0=e$ and $v_j=v$. Then for any $0\leq i\leq j$,
we have $v_{i-1}v_i^{-1}\in \{a_i^2\ |\ 1\leq i\leq k\}$.
Noting that the elements in $ \{a_i^2\ |\ 1\leq i\leq k\}$ are all involutions and commute with each other,
it follows that $v=a_{i_1}^2a_{i_2}^2\ldots a_{i_j}^2$
for some $1\leq i_1,\ldots, i_j\leq k$. Consequently,
$v\in \lg a_i^2\ |\ 1\leq i\leq k\rg$. Set $H=\lg a_i^2\ |\ 1\leq i\leq k\rg$.
The above argument gives that $V(A)\subseteq H$. Clearly, $H\cong\mz_2^k$,
and the subgraph of $\G_{k,\ell}$ induced by $H$ is just the Cayley graph $\Cay(H,\{a_i^2\ |\ 1\leq i\leq k\})$,
which is isomorphic to $Q_k$.

For any two distinct vertices $u,v\in V(A)$, if $u,v$ have a common neighbor, say $w$, in
$V(\G_{k,\ell})\setminus H$, then $w=s_1u=s_2v$ for some $s_1,s_2\in \Omega_{k,\ell}\setminus\{a_i^2\ |\ 1\leq i\leq k\}$,
and hence $s_1s_2=uv^{-1}\in \lg a_i^2\ |\ 1\leq i\leq k\rg$.
However, it is easy to check that the product of any two distinct elements in
$\Omega_{k,\ell}\setminus\{a_i^2\ |\ 1\leq i\leq k\}$ can not be in $\lg a_i^2\ |\ 1\leq i\leq k\rg$.
Thus, $s_1=s_2$. Since $w=s_1u=s_2v$, one has $u=v$, a contradiction.
Thus, any two distinct vertices of $A$ do not share a common neighbor in $V(\G_n)\setminus H$.
It follows that $|N_{V(\G_{k,\ell})\setminus H}(A)|=|A|(n-k)$,
and so $$|N_{H}(A)|=f_n(g)-|A|(n-k)=f_k(g).$$
By Lemma~\ref{star}, $Q_k[A]=\G_{k,\ell}[A]$ is a star, as required.\medskip

\f{\bf Case~2}\ There is an edge $\{x,y\}$ of $A$ such that $xy^{-1}\notin\{a_i^2\ |\ 1\leq i\leq k\}$.

Let $xy^{-1}=b$. Then $b\in\Omega_{k,\ell}$ and $x=by$. Set $H=\lg\Omega_{k,\ell}-b\rg$, $\G_{k,\ell}^{0}=\Cay(H, \Omega_{k,\ell}-b)$
and $\G_{k,\ell}^1=\G_{k,\ell}[Hb]$.
It is easy to see that $H\cong G_{k-1,\ell+2}$ or $G_{k,\ell-1}$, and
$\G_{k,\ell}^{0}\in \mathcal{G}_{n-1}$.
In view of the fact that
the map $R(b): v\mapsto vb, \forall v\in G_{k,\ell}$ is an automorphism of $\G_{k,\ell}$, it follows that
$\G_{k,\ell}^1=\G_{k,\ell}[Hb]\cong \G_{k,\ell}^0$. Also, since $b\notin\{a_i^2\ |\ 1\leq i\leq k\}$,
$H$ has index $2$ in $G_{k,\ell}$, and so $H$ has two cosets, namely, $H$ and $Hb$, in $G_{k,\ell}$.
Without loss of generality, assume that $y\in H=V(\G_{k,\ell}^0)$. Then $x=by\in bH=Hb=V(\G_{k,\ell}^1)$.
So $U_i=V(\G_{k,\ell}^i)\cap V(A)$ is non-empty for $i=0,1$. Without loss of generality, assume that $|U_0|\leq |U_1|$.
Set $|U_0|=N$. Then $N\leq \lfloor\frac{g+1}{2}\rfloor$.
By Proposition~\ref{basic}, $|N_{\G_{k,\ell}^0}(U_0)|\geq f_{n-1}(N-1)$
and $|N_{\G_{k,\ell}^1}(U_1)|\geq f_{n-1}(g-N)$.
Note that $N_{\G_{k,\ell}^0}(U_0)\cup N_{\G_{k,\ell}^1}(U_1)\subseteq N_{\G_{k,\ell}}(V(A))$ and $N_{\G_{k,\ell}^0}(U_0)\cap N_{\G_{k,\ell}^1}(U_1)=\emptyset$.
It follows that $|N_{\G_{k,\ell}^0}(U_0)|+|N_{\G_{k,\ell}^1}(U_1)|\leq |N_{\G_{k,\ell}}(V(A))|=f_n(g)$,
and hence
\begin{equation}\label{eq-1}
\begin{array}{lll}
0&\geq &|N_{\G_n^0}(U_0)|+|N_{\G_n^1}(U_1)|-|N_{\G_n}(V(A))|\\
&\geq & f_{n-1}(N-1)+f_{n-1}(g-N)-f_n(g)\\
&=&-(N-1)(N-g).
\end{array}
\end{equation}
Remember that $1\leq N\leq \lfloor\frac{g+1}{2}\rfloor$. If $N>1$,
then $N<g$ and so $0\geq -(N-1)(N-g)>0$, a contradiction. Thus, $N=1$ and so $|U_1|=g$.
Let $U_0=\{v\}$. Since $v$ has only one neighbor in $V_1$,
$v$ is a pendant vertex of $A$, and so $\G_{k,\ell}[U_1]$ is connected.
Again, since $N=1$, by the above equation (\ref{eq-1}),
we have $|N_{\G_{k,\ell}^1}(U_1)|=f_{n-1}(g-1)$ and $|N_{\G_{k,\ell}}(V(A))|=|N_{\G_{k,\ell}^0}(U_0)|+|N_{\G_{k,\ell}^1}(U_1)|$.
Then $|N_{\G_{k,\ell}}(U_1)|=|N_{\G_{k,\ell}^1}(U_1)|+|N_{\G_{k,\ell}^0}(U_1)|=f_{n}(g-1)$.
By the induction hypothesis,
we have $\G_{k,\ell}^1[U_1]$ is a star or a $3$-path.
The equation $|N_{\G_{k,\ell}}(V(A))|=|N_{\G_{k,\ell}^0}(U_0)|+|N_{\G_{k,\ell}^1}(U_1)|$ implies that
$N_{\G_{k,\ell}^0}(U_1-\a)\subseteq N_{\G_{k,\ell}^0}(v)$, where $\a$ is the neighbor of $v$ in $\G_{k,\ell}^1$.

Suppose that $\G_{k,\ell}[U_1]$ is a $3$-path $u_1u_2u_3u_{4}$.
Assume that the neighbor of $u_i$ in $\G_{k,\ell}^0$ is $w_i$ for each $1\leq i \leq 4$.
Then $v\in\{w_i\ |\ 1\leq i\leq 4\}$.

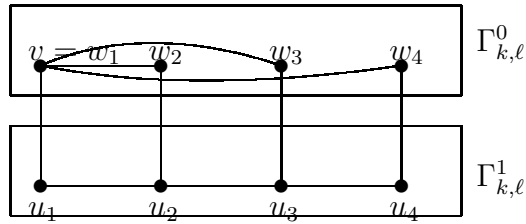
\begin{figure}[ht]
\begin{center}
\unitlength 4mm
\begin{picture}(15,8)

\put(0,4){\line(1, 0){15}}
\put(0,4){\line(0, 1){3}}
\put(15,7){\line(-1, 0){15}}
\put(15,7){\line(0, -1){3}}
\put(15.5,5.5){$\G_{k,\ell}^0$}

\put(0,0){\line(1, 0){15}}
\put(0,0){\line(0, 1){3}}
\put(15,3){\line(-1, 0){15}}
\put(15,3){\line(0, -1){3}}
\put(15.5,1){$\G_{k,\ell}^1$}

\put(1, 1){\circle*{0.4}} \put(0.6, 0){$u_1$}
\put(5, 1){\circle*{0.4}} \put(4.6, 0){$u_2$}
\put(9, 1){\circle*{0.4}} \put(8.6, 0){$u_3$}
\put(13, 1){\circle*{0.4}} \put(12.6, 0){$u_4$}

\put(1,1){\line(1, 0){12}}

\put(1,1){\line(0, 1){4}}\put(5,1){\line(0, 1){4}}
\put(9,1){\line(0, 1){4}}\put(13,1){\line(0, 1){4}}
\put(1, 5){\circle*{0.4}}\put(5, 5){\circle*{0.4}}
\put(9, 5){\circle*{0.4}}\put(13, 5){\circle*{0.4}}
 \put(0.6, 5.2){$v=w_1$} \put(4.6, 5.2){$w_2$}
 \put(8.6, 5.2){$w_3$}\put(12.6, 5.2){$w_4$}

\put(1,5){\line(1, 0){4}}

 \qbezier(1, 5)(4.5, 6.5)(9, 5)
  \qbezier(1, 5)(6, 4)(13, 5)

\end{picture}
\end{center}\vspace{-.2cm}
\caption{Illustration of the proof for the case $v=w_1$} \label{Fig-4.1}
\end{figure}

Let $v=w_1$ or $w_4$. Without loss of generality, assume that $v=w_1$(see Figure~\ref{Fig-4.1}).
Suppose that $u_1$ and $w_3$
has a common neighbor, say $x$, such that $x\neq v$. Since the edges between $\G_{k,\ell}^0$ and $\G_{k,\ell}^1$
are vertex-disjoint, we must have $x=u_3$, forcing $(u_1,u_2,u_3)$ is a triangle, contrary to the fact that $\G_{k,\ell}$
has girth $4$ (Proposition~\ref{com-zhu}). Thus, $v$ is the unique common neighbor of $u_1$ and $w_3$.
By Lemma~\ref{common neighbor}, $\{u_1,w_3\}=\{b_iv,a_ib_iv\}$ for some $1\leq i\leq s$.
Clearly, $v$ is also a common neighbor of $u_1$ and $w_4$.
By  Lemma~\ref{common neighbor}, if $u_1$ and $w_4$ have only one common neighbor,
then we also have $\{u_1,w_4\}=\{b_iv,a_ib_iv\}$, which is impossible. Thus, $u_1$ and $w_4$
must have two common neighbors,
and so $u_4$ must be a common neighbor of $u_1$ and $w_4$.
Thus, $(u_1,u_2,u_3,u_4)$ is a $4$-cycle. By Proposition~\ref{com-zhu}, any two vertices of $\G_{k,\ell}$ have at most two common neighbors.
It follows that
$|N_{\G_{k,\ell}^1}(U_1)|=4(n-3)>f_{n-1}(3)=4n-13$, a contradiction.

\begin{figure}[ht]
\begin{center}
\unitlength 4mm
\begin{picture}(18,8)

\put(0,4){\line(1, 0){17}}
\put(0,4){\line(0, 1){3}}
\put(17,7){\line(-1, 0){17}}
\put(17,7){\line(0, -1){3}}
\put(17.5,5.5){$\G_{k,\ell}^0$}

\put(0,0){\line(1, 0){17}}
\put(0,0){\line(0, 1){3}}
\put(17,3){\line(-1, 0){17}}
\put(17,3){\line(0, -1){3}}
\put(17.5,1){$\G_{k,\ell}^1$}

\put(1, 1){\circle*{0.4}} \put(0, 0){{\footnotesize$u_1=bg_1v$}}
\put(5, 1){\circle*{0.4}} \put(4.2, 0){{\footnotesize$u_2=bv$}}
\put(9, 1){\circle*{0.4}} \put(7.8, 0){{\footnotesize$u_3=b'bb'v$}}
\put(13, 1){\circle*{0.4}} \put(12.6, 0){{\footnotesize$u_4=bb'v$}}

\put(1,1){\line(1, 0){12}}

\put(1,1){\line(0, 1){4}}\put(5,1){\line(0, 1){4}}
\put(9,1){\line(0, 1){4}}\put(13,1){\line(0, 1){4}}
\put(1, 5){\circle*{0.4}}\put(5, 5){\circle*{0.4}}
\put(9, 5){\circle*{0.4}}\put(13, 5){\circle*{0.4}}
 \put(0, 5.2){{\footnotesize$w_1=g_1v$}} \put(4.0, 5.2){{\footnotesize$w_2=v$}}
 \put(8, 5.2){{\footnotesize$w_3=a_i^2v$}}\put(12.6, 5.2){{\footnotesize$w_4=b'v$}}

\put(1,5){\line(1, 0){4}}
\put(5,5){\line(1, 0){4}}

  \qbezier(5, 5)(9, 4)(13, 5)

\end{picture}
\end{center}\vspace{-.2cm}
\caption{Illustration of the proof for the case $v=w_2$} \label{Fig-4.2}
\end{figure}
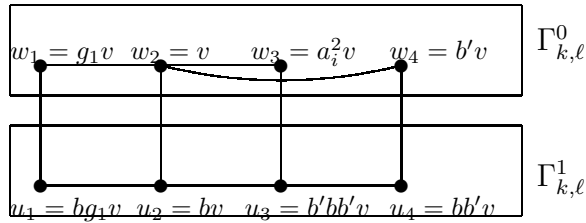

Let $v=w_2$ or $w_3$. Without loss of generality, assume that $v=w_2$(see Figure~\ref{Fig-4.2}).
Suppose that $w_4$ and $u_2$ have a common neighbor, say $w$, different from $v$.
If $w\in V(\G_{k,\ell}^0)$, then since $u_2$ has a unique neighbor in $V(\G_{k,\ell}^0)$,
we must have $w=v$, a contradiction. If $w\in V(\G_{k,\ell}^1)$, then since $w_4$ has a unique neighbor in $V(\G_{k,\ell}^1)$,
we must have $w=u_4$. It follows that $(u_2,u_3,u_4=w)$ is
a triangle, contrary to the fact that $\G_{k,\ell}$ has girth $4$ (Proposition~\ref{com-zhu}).
Thus, $v$ is the unique common neighbor of $w_4$ and
$u_2$. With a similar argument, we have $u_4$ is the unique
common neighbor of $u_3$ and $w_4$, and $w_4$ is the unique
common neighbor of $w_2$ and $u_4$.

By Lemma~\ref{common neighbor}, $\{w_4,u_2\}=\{b_iv,a_ib_iv\}$ for some $1\leq i\leq k$.
Let $w_4=b'v$ for some $b'\in\Omega_{k,\ell}$. It is easy to see that
for each $u\in H=V(\G_{k,\ell}^0)$, $bu$ is the unique neighbor of $u$
in $Hb=V(\G_{k,\ell}^1)$. Thus, $u_2=bw_2=bv$, and so $\{b'v,bv\}=\{w_4,u_2\}=\{b_iv,a_ib_iv\}$.
It follows that $\{b,b'\}=\{b_i,a_ib_i\}$.
Also, $u_4=bw_4=bb'v$. As $u_4$ is the unique
common neighbor of $u_3$ and $w_4$, by  Lemma~\ref{common neighbor},
$\{w_4,u_3\}=\{bu_4,b'u_4\}=\{b'v,b'bb'v\}$ because $w_4=bu_4$.
Consequently, $u_3=b'u_4=b'bb'v$, and so $w_3=bu_3=(bb')^2v=a_i^2v$ (note that $bb'=a_i$ or $a_i^{-1}$).
Set $w_1=g_1v$ with $g_1\in\Omega_{k,\ell}$. Then $u_1=bg_1v$.
Note that $\{w_3,w_4,u_2\}=\{a_i^2v, a_ib_iv,b_iv\}$.
So, $g_1\notin\{a_i^2,a_ib_i,b_i\}$. Since $g_1\in\Omega_{k,\ell}$, one has $g_1\notin\lg a_i,b_i\rg$, and so $g_1$ commutes
with the elements in $\lg a_i,b_i\rg$. In particular,
$g_1$ commutes with $b$ and $b'$.
Since $u_2=bv$ is adjacent to $u_3=b'bb'v$, there exists
an $s\in\Omega_{k,\ell}$ such that $su_2=u_3$, namely, $sb=b'bb'$.
It follows that $s=(b'b)^2=a_i^2$. Now one may see that
$u_1(=g_1u_2)$ and $u_3(=su_2)$ have a common neighbor $g_1su_2$ which is different from $u_2$.
Moreover, $u_2(=a_i^2u_3)$ and $u_4(=b'u_3)$ also have a common neighbor $a_i^2b'u_3$ which is different from $u_3$.
By Proposition~\ref{com-zhu}, any two vertices of $\G_{k,\ell}$ have at most two common neighbors.
Since $|N_{\G_{k,\ell}^1}(U_1)|=f_{n-1}(3)=4n-13$, $u_1$ and $u_4$ must have a unique common neighbor, say $w$.
Again by Lemma~\ref{common neighbor}, we have $\{u_1,u_4\}=\{b_jw,a_jb_jw\}$, and so
$u_1u_4^{-1}\in\{a_j,a_j^{-1}\}$ for some $1\leq j\leq k$.
This implies that $u_1u_4^{-1}$ has order $4$. However, $u_1u_4^{-1}=bg_1v(bb'v)^{-1}=g_1a_i^2b'$ is an involution,
a contradiction.

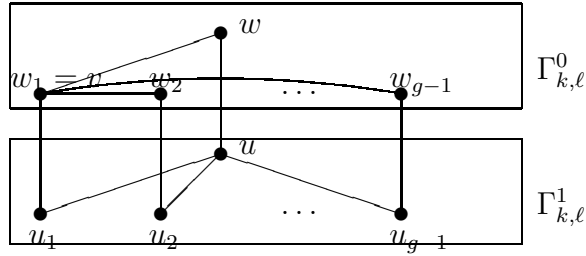
\begin{figure}[ht]
\begin{center}
\unitlength 4mm
\begin{picture}(18,10)

\put(0,4.5){\line(1, 0){17}}
\put(0,4.5){\line(0, 1){3.5}}
\put(17,8){\line(-1, 0){17}}
\put(17,8){\line(0, -1){3.5}}
\put(17.5,5.5){$\G_{k,\ell}^0$}

\put(0,0){\line(1, 0){17}}
\put(0,0){\line(0, 1){3.5}}
\put(17,3.5){\line(-1, 0){17}}
\put(17,3.5){\line(0, -1){3.5}}
\put(17.5,1){$\G_{k,\ell}^1$}

\put(1, 1){\circle*{0.4}} \put(0.6, 0){{$u_1$}}
\put(5, 1){\circle*{0.4}} \put(4.6, 0){{$u_2$}}
\put(9, 1){$\ldots$} 
\put(13, 1){\circle*{0.4}} \put(12.6, 0){{$u_{g-1}$}}
\put(7, 3){\circle*{0.4}}\put(7.6, 3){{$u$}}

\put(7,3){\line(-1, -1){2}}
\put(7,3){\line(-3, -1){6}}
\put(7,3){\line(3, -1){6}}
\put(7,3){\line(0, 1){4}}

\put(1,1){\line(0, 1){4}}\put(5,1){\line(0, 1){4}}
\put(13,1){\line(0, 1){4}}
\put(1, 5){\circle*{0.4}}\put(5, 5){\circle*{0.4}}
\put(9, 5){$\ldots$}\put(13, 5){\circle*{0.4}}
 \put(0, 5.2){{$w_1=v$}} \put(4.6, 5.2){{$w_2$}}
 \put(12.6, 5.2){{$w_{g-1}$}}
 \put(7, 7){\circle*{0.4}}\put(7.6, 7){{$w$}}

\put(1,5){\line(1, 0){4}}
\put(1,5){\line(3, 1){6}}

  \qbezier(1, 5)(7, 6)(13, 5)

\end{picture}
\end{center}\vspace{-.2cm}
\caption{Illustration of the proof for the case $v=w_1$} \label{Fig-4.3}
\end{figure}

Now we assume that $\G_{k,\ell}[U_1]$ is a star with vertex-set $U_1=\{u,u_1,\ldots,u_{g-1}\}$
and edge-set $E(\G_{k,\ell}[U_1])=\{\{u,u_i\}\ |\ 1\leq i\leq g-1\}$.
If $v$ is adjacent to $u$,
then $A$ is a star. Assume that $v$ is adjacent to some $u_i$.
Without loss of generality, let $i=1$.
Assume that the neighbor of $u_i$ in $\G_{k,\ell}^0$ is $w_i$ for each $1\leq i \leq g-1$,
and the neighbor of $u$ in $\G_{k,\ell}^0$ is $w$.
Then $v=w_1$.
If $g\leq 2$, then $A$ is also a star.
If $g=3$, then $A$ is a $3$-path. Suppose $g>3$(see Figure~\ref{Fig-4.3}).
Since $N_{\G_{k,\ell}^0}(U_1-u_1)\subseteq N_{\G_{k,\ell}^0}(v)$, for any $u_j\in U_1-\{u,u_1\}$,
the neighbor $w_j$ of $u_j$ in $\G_{k,\ell}^0$ is also adjacent to $v$, and
hence $(v,u_1,u,u_j,w_j)$ is a cycle of length $5$.
Suppose that $v$ is not the unique common neighbor
of $u_1$ and $w_j$ for some $2\leq j\leq g-1$. Then $u_1$ and $w_j$ have another
common neighbor say $x$, different from $v$. Note that the edges between $V(\G_n^0)$ and
$V(\G_n^1)$ form a perfect matching. If $x\in V(\G_{k,\ell}^0)$, then $x=v$, a contradiction.
If $x\in V(\G_{k,\ell}^1)$, then $x=u_j$, and so $(u_1,u,u_j)$ is a triangle, contrary to the
fact that $\G_{k,\ell}$ has girth $4$ (Proposition~\ref{com-zhu}). Thus, $v$ is the unique common neighbor
of $u_1$ and $w_j$ for each $2\leq j\leq g-1$. By Lemma~\ref{common neighbor},
we have $\{u_1,w_j\}=\{b_{i_j}v,a_{i_j}b_{i_j}v\}$ for some $0\leq i_j\leq k$. Since $g>3$,
one has $g-1>2$, and hence
$\{u_1,w_2\}=\{b_{i_2}v,a_{i_2}b_{i_2}v\}$ and $\{u_1,w_{g-1}\}=\{b_{i_{g-1}}v,a_{i_{g-1}}b_{i_{g-1}}v\}$.
Since $u_1\in\{b_{i_{g-1}}v,a_{i_{g-1}}b_{i_{g-1}}v\}\cap\{b_{i_2}v,a_{i_2}b_{i_2}v\}$, by the structure of the group
$G_n$ we must have $i_2=i_{g-1}$.
It follows that $\{u_1,w_2\}=\{u_1,w_{g-1}\}$, and hence $w_2=w_{g-1}$, a contradiction.
\hfill\qed

\begin{theorem}\label{cor}
If $2k+1\leq g\leq 3k-4$, then $\k_g(\G_{k,0})>f_{3k}(g)$.
\end{theorem}\medskip

\demo By Theorem~\ref{extra-0}, $\k_g(\G_{k,0})\geq f_{3k}(g)$.
Suppose on the contrary that $\k_g(\G_{k,0})=f_{3k}(g)$. Let $S$ be
a minimum $R_g$-cutset of $\G_{k,0}$. By~Theorem~\ref{extra-0},
$\G_{k,0}-S$ has two components, one of which, say $A$, has $g+1$ vertices.
By the minimality of $S$, we see that $S=N_{\G_{k,0}}(V(A))$. Since
$|S|=\k_g(\G_{k,0})=f_{3k}(g)$, by Lemma~\ref{component},
$A$ is isomorphic to $K_{1,g}$ or $P_3$. Since $2k+1\leq g\leq 3k-4$, one has
$k\geq 5$, and so $g+1\geq 12$. It follows that $A\cong K_{1,g}$. By the
vertex-transitivity of $\G_{k,0}$, we may assume that the identity element $e$ of
$G_{k,0}$ is a vertex of $A$, and $V(A)=\{e,s_1,s_2,\ldots,s_g\}$.
By Proposition~\ref{com-zhu}, any two vertices of $\G_{k,0}$ have at most two common neighbors.
Since $|N_{\G_{k,0}}(V(A))|=f_{3k}(g)$, by an easy calculation,
we see that for any $s_i,s_j\in V(A)-\{e\}$, they have exactly two common neighbors.
Clearly, $V(A)-\{e\}\subseteq\Omega_{k,0}$.
Since $g\geq 2k+1$, there must exist $1\leq i\leq k$ such that
$a_ib_i,b_i\in V(A)$. So, $a_ib_i$ and $b_i$ have two common neighbors.
Clearly, $e$ is a common neighbor of them. Let $u\neq e$ be another
common neighbor of $a_ib_i$ and $b_i$. Then $u=xa_ib_i=yb_i$ for some
$x,y\in\Omega_{k,0}$. It follows that $xy=a_i$.
In view of the fact that $\Omega_{k,0}=\{a_i^2,b_i,a_ib_i\ |\ 1\leq i\leq k\}$,
we must have $x=a_ib_i$ and $y=b_i$, forcing $u=e$, a contradiction.\hfill\qed

\f{\bf Remark~2}\ This theorem implies that \cite[Theorem~3.4]{Zhu}, which states that
the $g$-extra connectivity of an $n$-dimensional HL-network $G$ is $f_n(g)$ for $0\leq g\leq n-4$, is not true.
In fact, in the proof of \cite[Theorem~3.4]{Zhu}, the authors
first assumed that $G=G_0\oplus G_1$, and then picked a subgraph, say $A$, of $G$ such that
$A\cong K_{1,g}$ and $A\subseteq G_0$. Then they claimed that ``it is straightforward to see
that $|N_G(A)|=n(g+1)-\frac{g(g+3)}{2}$." However, from the proof of Theorem~\ref{cor} we see that this claim is not true.

\section{Conclusion}
As one of
the novel factors for measuring the reliability and fault tolerance of networks,
the $g$-extra connectivity of HL-networks have also been studied by some
authors. From \cite{YM} we see that $\k_g(Q_n)=f_n(g)=n(g+1)-\frac{1}{2}g(g+3)$
where $n\geq 4$ and $0\leq g\leq n-3$. In view of the fact that the hypercubes $Q_n$
belong to the class of HL-networks, an interesting problem is:
For an $n$-dimensional HL-network $G$, does $\k_g(G)=f_n(g)$ hold for $n\geq 4$ and $0\leq g\leq n-3$?

For this problem, some partial answers have been known. Several authors proved that
the answer is positive for the case when $g\leq 3$ (see \cite{FH,Zhu2008,XZX,CH}).
In this paper, a subclass of $n$-dimensional HL-networks with $g$-extra connectivity greater than $f_n(g)$
is presented, and so a negative answer to the above problem is given.

Moreover, we also prove that for $n\geq 5$ and $0\leq g\leq n-3$, if $n$-dimensional HL-network
$X_n$ has a $R_g$-cutset $S$ with $|S|\leq f_n(g)$, then $|S|=f_n(g)$ and $X_n$ is hyper-$\k_g$.
This enables us to obtain a lower bound on $\k_g(G)$, namely, $\k_g(X_n)\geq f_n(g)$.
In addition, we also give a sufficient condition for the lower bound being attainable.
Applying this, we first present a short proof for the main results of \cite{CH} which shows that
the answer to the above problem is positive for the case when $g=2$ or $3$, and then
determine the $g$-extra connectivity of the varietal hypercubes for some specific $g$.

In view of the above facts, an interesting problem is: Determine the
smallest $g$ such that $\k_g(G)=f_n(g)$ holds for all $n$-dimensional HL-networks $G$.
This is a topic for our future effort.

\section*{Acknowledgment}

This work was supported by the National Natural Science Foundation of China (11271012)  and the Fundamental
Research Funds for the Central Universities (2015JBM110).

\end{document}